\documentclass[a4paper,12pt]{amsart}

 \usepackage{graphicx}

\input xy
\xyoption{all}
\usepackage{amssymb}

\usepackage[left=2.2cm,top=2.5cm,right=2.2cm,bottom=2.5cm]{geometry}

\email{roysland@math.uio.no}
\newtheorem{theorem}{Theorem}[section]
\newtheorem{lemma}[theorem]{Lemma}
\newtheorem{proposition}[theorem]{Proposition}

\newtheorem*{acknowledgments}{Acknowledgments}
\theoremstyle{definition}
\newtheorem{definition}[theorem]{Definition}

\theoremstyle{remark}
\newtheorem{remark}[theorem]{Remark}
\numberwithin{equation}{section}

\def\R{\mathbb{R}}
\def\N{\mathbb{N}}

\def\T{\mathbb{T}}

\def\B{\mathcal{B}}

\def\P{\mathcal{P}}

\def\-{^{-1}}

\def\K{\mathcal{K}}

\def\C{\mathbb{C}}
\def\Z{\mathbb{Z}}

\begin{document}

\title{Symmetries in  Projective multiresolution analyses}

\author{Kjetil R{\o}ysland}
\footnote{Research supported in part by the Research Council of Norway, project number NFR 154077/420. Some of the final work was also done with the support from the project NFR 170620/V30.}
\address{[Kjetil R\o ysland]University of Oslo\\
Department of Mathematics\\
PO Box 1053, Blindern\\
NO-0316 Oslo\\
Norway}
\email{roysland@math.uio.no}


\begin{abstract}
We give an equivariant version of Packer and Rieffel's theorem on sufficient conditions for the existence of orthonormal wavelets in projective multiresolution analyses. The scaling functions that generate a projective multiresolution analysis are supposed to be invariant with respect to some finite group action. We give sufficient conditions for the existence of wavelets with similar invariance.
\end{abstract}
\maketitle
\thispagestyle{empty}

\section{Introduction} 
Let $G$ be a full rank lattice in $\R^n$,  i.e. a subgroup such that  $G \simeq \Z^n$, 
and let    $X = \R^n / G$. We let $p : \R^n \rightarrow X$ denote the quotient map. 
Following \cite{PR1},  we let $\Xi$ denote the set of $f \in C_b(\R^n)$ such that 
$\sum_{g \in G} |f |^2(x-g)$ defines  a continuous function on $X$. We define 
$\langle \cdot, \cdot \rangle' : \Xi \times \Xi \rightarrow C(X )$
by
$$
\langle \zeta, \eta \rangle' \circ p(x) = \sum_{g \in G} \zeta(x-g) \overline \eta(x - g). 
$$
As shown in \cite{PR1}, $\Xi$ and $ \langle \cdot, \cdot \rangle'$ form a  $C(X)$-Hilbert module. 

 Let $A \in GL_n(\R)$ such that 
$AG \subset G$ and if $\lambda$ is an eigenvalue of $A$ then $|\lambda | > 1$. 
We let $q = | \det A |$ and define $U \in B(\Xi)$ such that 
$U \zeta = q^{-1/2} \zeta \circ A\-$. 
A projective multiresolution analysis is a family of 
subspaces  $\{V_j\}_{j \in \Z}$  of $\Xi$ such that
\begin{enumerate}
  \item $V_j \subset V_{j+1}$,
  \item $U V_j =  V_{j+1}$,
  \item $V_0$ is a projective submodule in $\Xi$,
  \item $\cup_j V_j$ is dense in $\Xi$,
  \item $\cap_j V_j = \{0\}$.
\end{enumerate}
These conditions are redundant in the sense that we get a projective multiresolution analysis if we have 
a projective submodule $V \subset \Xi$, 
such that $V \subset U V$ and $\cup_j U^j V$ is dense in $\Xi$.
Projective multiresolution analyses are closely related to multiresolution analyses. 
If the core space is free, this relation is particularly simple.

\begin{proposition}
Suppose  $\{V_j\}_{j \in \Z}$ form a projective multiresolution analysis.  
Moreover, let  $\phi_1, \dots, \phi_d$ be  an orthonormal $C(X)$-basis for $V_0$. 
Let $G^\bot$ denote  the dual lattice of $G$, i.e. $G^\bot = \{ y \in \R^n | \langle x , y \rangle \in \Z 
\text{ for every } x \in G\}$.
 The elements of the set  $\{ \check \phi_j( \cdot - g) | 1 \leq j \leq d, ~ g \in G^\bot \} $ form an orthonormal family in $L^2(\R^n)$. 
Let $W_0$ denote the minimal closed space that contains this family, define  
$\check U \in B(L^2(\R^n))$ as $\check U f = \sqrt q  f \circ A^{T}$ and let $W_j = \check U^j W_0$. 

The triple  $(\check U,\{W_j\}_{j \in \Z}), \{\check \phi_i  \}_i$ defines  a multiresolution analysis in the sense that
\begin{enumerate}
\item $\check U W_j =  W_{j+1}$,
\item $W_j \subset W_{j+1}$,
\item $\cap_{j \in \Z} W_j = \{0\}$,
\item $\cup_{j \in \Z} W_j$ is dense in $L^2(\R^n)$,
\item The elements in $\{ \check \phi_j( \cdot - g) | 1 \leq j \leq d, ~ g \in G^\bot \} $ form an orthonormal basis for $W_0$.
\end{enumerate}    

\begin{proof}

Let $\mu$ denote the normalized Haar measure on $X$ and let $\tilde \mu$ denote the 
Haar measure on $\R^n$ such that 
$
\int_X \sum_{g \in G} f ( x - g) d\mu(x) = \int_{\R^n} f d \tilde \mu$
for every $f \in C_c(\R^n)$. 

If $\zeta , \eta  \in \Xi$, we see that $\int_{\R^n} f \overline g d \tilde \mu = \int_X \langle f,g \rangle' d\mu $, 
so $\Xi \subset L^2(\R^n , \tilde \mu)$.
 The dual lattice $G^\bot$,
is a full rank lattice in $\R^n$. 
We define a pairing 
 $ ( \cdot ,\cdot ) : X \times G ^\bot \rightarrow \T$ as  
$(x,g) = e^{2 \pi i \langle y, g \rangle}$ for any $y$ in the equivalence class of $x$. 
This pairing yields an isomorphism  $\widehat{\R^n /G} \simeq G^\bot$, 
see \cite[Theorem 4.39]{Folland}.

If  $\zeta_1, \dots, \zeta_r \in  \Xi$ and $g \in G^\bot$ then  
\begin{align*}
  \int_{\R^n} \check \zeta_i( y - g)  \overline{ \check   \zeta}_j(y ) d \tilde \mu(y)  
  & = \int_{\R^n} \zeta_i(y) \overline{\zeta}_j ( y) e^{-i 2 \pi \langle y, g \rangle    }     d \tilde \mu (y) = \int_X \langle \zeta_i , \zeta_j \rangle'(x)  \overline{(x,g)} d \mu(x)
\end{align*}
This implies that $\{ \check \zeta_j ( \cdot  - g) | 1 \leq j \leq r, g \in G^\bot \}\subset L^2(\R^n, \tilde \mu)$
 are orthonormal 
if and only if $\langle \zeta_i, \zeta_j \rangle' = \delta_{i,j}$.

We see that if $V_0$ is finitely generated and free with an orthonormal $C(X)$-basis $\phi_1, \dots , \phi_d \in V_0$,
then the inverse Fourier transform yields  an ordinary multiresolution analysis in $L^2( \R^n , \tilde\mu)$, 
$( \{W_j\}_{j \in \Z},   \check U)$, such that $ W_{j+1} = \check UW_{j}$  and   
$\check \phi_1, \dots, \check \phi_d \in W_0$ are  scaling functions with mutually orthonormal $G^\bot$-translates.

\end{proof}
\end{proposition}

An orthonormal MRA wavelet family with respect to  $(\check U,\{W_j\}_{j \in \Z})$
is a family \\ $\check \psi_1, \dots, \check \psi_r \in W_1$ such that 
$\{ \check \psi_j ( \cdot - g) | 1 \leq j \leq r, ~ g \in G^\bot\}$
form an orthonormal basis for $W_1 \ominus W_0$.
As noted in \cite{PR2}, there will always exist an orthonormal MRA  wavelet family if 
the scaling functions of the MRA 
have orthonormal translates. Suppose in addition that
$\phi_1, \dots , \phi_d  \in \Xi$. 
Do there exist MRA wavelets such that $ \psi_j \in \Xi$ ? Packer and Rieffel 
gave  an answer to this question in  \cite{PR2}. 
 The starting point of their argument is the following observation about the quotient $C(X)$-module $(UV_0)/V_0$:

\begin{proposition}  \label{ortbase}
If $V_0$ is a free $C(X)$-module of  rank $d$, then 
$U V_0$ is a free $C(X)$-module of rank $dq$. We have the following isomorphism 
of $C(X)$-modules
$$(U V_0) /V_0   \oplus V_0 \simeq U V_0.$$
Moreover, $(U V_0) / V_0$ is free if and only if  there exists an orthonormal MRA wavelet family 
$\check \psi_1, \dots, \check \psi_{d(q-1)} \in W_1$ such that  
$\psi_1, \dots ,  \psi_{d(q-1)} \in UV_0$.

\begin{proof}
$A^T G^\bot$ is a subgroup in $G^\bot$ since $(A\- G)^\bot = A^T G ^\bot$. 
%
We define a map $( \cdot, \cdot ) : G/AG \times G^\bot  / A^T G^\bot \rightarrow \T$,  by 
$(x,y) = e^{2 \pi i \langle A\- x' , y' \rangle}$ where 
$x' \in G$, $[x']_{AG} = x$, $y' \in G^\bot $ and $[y']_{A^T G^\bot} = y$. First, we note that if 
$x_0 \in AG$ and $y_0 \in A^T G^\bot$, then $e^{2\pi i \langle A\- x' +A\-  x_0   , y' + y_0  \rangle} = 
e^{2 \pi i \langle A\- x' , y' \rangle}$, so the map $( \cdot, \cdot )$ is well defined. 
Moreover, we see that if $y_1 , y_2 \in G^\bot$ such that $y_1 - y_2 \notin A^T G^\bot$, 
there exists an $x \in G$ such that $e^{2 \pi i \langle A\- x, 
  y_1 \rangle} \neq e^{2 \pi i \langle A\- x, y_2 \rangle}$. 
This means that representatives for different cosets of $A^T G^\bot$ in $G^\bot$ give different 
characters on $G/AG$. This gives us an injective group homomorphism $G^\bot/A^T G^\bot  \rightarrow \widehat{G/AG}$. 
Finally we note that 
$| G^\bot  / A^TG^\bot | = |\det A^T| = |\det A |= | G/AG | = | \widehat {G/AG}|$ 
and our group homomorphism must be surjective. 
 
Let $ y_1, \dots, y_q$ be a system of representatives for $A^TG^\bot$ in $G^\bot$
and define $f_j = e^{2 \pi i \langle  A\- \cdot , y_j \rangle } \in C_b (\R^n)$, for $j = 1, \dots, q$.

Let $\B$ denote the algebra of bounded and continuous $AG$-periodic 
functions on $\R^n$. Now $C(X) \subset \B$ is a subalgebra with a conditional expectation
$P : \B \rightarrow C(X)$ defined as $Pf (p(x)) = q\- \sum_{j = 1    }^q   f ( x - A\-x_j)$
where $x_1, \dots , x_q \in G$ form a system 
of representatives for the cosets of $AG$ in $G$. 
A direct computation gives that $P(f_k \overline f_l)(p(x)) = f_k \overline f_l q\- \sum_j ([ x_j], [ y_k - y_l]     )$. The orthogonality of characters on compact groups implies  that $f_1, \dots, f_q$ form an orthonormal 
$C(X)$-basis for $\B$ with respect to 
the $C(X)$-module inner product $f, g \mapsto P( f \overline g)$. 

Whenever  $\zeta, \eta \in \Xi$, we have 
\begin{align*}
P (  \langle U\-  \zeta, U\-  \eta \rangle' \circ p A\-  )(x)  = & 
\sum_j \sum_{g \in G} (\zeta \overline \eta) ( x - x_j  - Ag) =   
\langle \zeta, \eta \rangle' \circ p(x)
\end{align*}
This implies that if 
$\phi_1, \dots, \phi_d \in V_0$ form an orthonormal $C(X)$-basis for $V_0$, 
then $\{f_k U \phi_l |1 \leq k \leq q,      1 \leq l \leq d\}$ form an orthonormal $C(X)$-basis for  
the module $U V_0$.

We let 
$Q : \zeta \mapsto \sum_{j = 1}^k  \langle \zeta, \phi_j \rangle \phi_j$,  the 
 orthogonal projection from $\Xi$, onto $V_0$. 
If $V_0 \subset U V_0$, then $U V_0 \simeq V_0 \oplus ( 1-Q) U V_0$ and 
$U V_0/V_0 \simeq  (1-Q)U V_0$ form a projective and stably free $C(X)$-module.

An orthonormal MRA wavelet family 
$\check \psi_1, \dots \check \psi_r$ such that $\psi_1, \dots , \psi_r \in \Xi$ 
would form an orthonormal $C(X)$-basis for $UV_0 \ominus V_0 \simeq UV_0 /V_0$, so  $r = d(q-1)$.
Conversely, suppose  $U V_0/V_0$ is free and  let $H : \oplus^{d(q-1)} C(X) \rightarrow (1-Q)U V_0$ 
be a $C(X)$-linear isomorphism. 
We equip $\oplus^{d(q-1)} C(X)$ with the ordinary $C(X)$-valued inner product. Now   
$H$ is adjointable with adjoint $H^* \zeta = \sum_j \langle \zeta, H e_j  \rangle_{\oplus^{d(q-1)} 
  C(X)}e_j $ where $\{e_j\}$ are  the standard orthonormal vectors. Moreover, $H$
is a bounded and  surjective operator between Banach spaces, so by the open mapping theorem, 
we obtain a $\delta >0$ such that for every $\zeta \in (1-Q)U V_0$ there exists an $\eta \in \oplus^{d(q-1)} C(X)$ such that
$H \eta = \zeta $ and $\| \eta \| \leq \delta\- \|\zeta \|$. Following \cite[3.2]{Lance} we see that 
\begin{align*}
  \|\zeta \|^2 & = \| \langle H \eta ,  \zeta \rangle' \| = 
  \| \langle \eta , H^* \zeta \rangle_{\oplus ^{d(q-1)  } C(X)} \| 
  \leq \| \eta \| \| H^* \zeta \|   
   =    \| \eta \| \| \langle \zeta, HH^* \zeta \rangle' \|^{1/2} \\ \leq  & \| \eta \| \| \zeta \|^{1/2} 
  \| HH^* \zeta \| ^{1/2}  = \delta\- \| \zeta \|^{3/2} \| HH^* \zeta \|^{1/2}
\end{align*}
i.e. $\| \zeta \| \leq \delta^{-2} \| HH^* \zeta \|$ for every $\zeta \in (1-Q)U V_0E$. Proposition \cite[3.1]{Lance}
states  that this condition implies that $HH^*$ is invertible, so we can apply spectral calculus for selfadjoint operators
to obtain a new operator $(HH^*)^{-1/2}H : \oplus^{ d(q-1) }C(X) \rightarrow (1-Q)U V_0$. 
Now $\psi_j = (HH^*)^{-1/2} H e_j $ form an orthonormal $C(X)$-basis for $(1-Q)U V_0$, i.e. we obtain an 
orthonormal MRA wavelet family $\check \psi_1, \dots , \check \psi_{d(q-1)}$ with $\psi_1, \dots, \psi_{d(q-1)} \in \Xi$.
\end{proof}
\end{proposition}

 We have that the quotient $(U V_0 )/V_0$ is a projective $C(X)$-module. 
The Serre-Swan theorem states that the complex vector bundles over $X$ and 
the finitely generated projective $C(X)$-modules are equivalent as categories. Moreover, every finitely generated and projective $C(X)$-module can be realized as the sections in a complex vector bundle over $X$. A free module now corresponds to the sections in a trivial bundle.   

Packer and Rieffel \cite{PR2} used this and  a cancellation result for vector bundles to 
say that the quotient  $(U V_0)/ V_0$ is free,
and hence allows  orthonormal MRA wavelets in $\Xi$ when 
$n < 2 q -1   $.

We follow their ideas, but in addition we will also consider some actions of  finite groups 
on  the projective multiresolution analysis. 
We  give an  equivariant version of the  cancellation result. 
This is used to  give sufficient conditions for the existence of symmetric MRA wavelets in $\Xi$.
In the last section, we classify these actions when $n = 2$.

\section{Symmetries}
\begin{definition} \label{affdef}
We will say that a non trivial and finite subgroup $H \subset GL_n(\Z)$ is affiliated to $A$ if it satisfies the following properties:  
\begin{enumerate}
\item $h  G = G$, \label{well}
\item $h A = A h$, \label{comm}
\item $( h-1) G \subset AG$, \label{triv}
\end{enumerate}
for every $h \in H$.
\end{definition}
If $n = 1$, then  $A = \pm 2$ and $H = \{\pm 1\}$ are the only examples. 
In the last section, we will give a list of all the possible pairs,  up to similarity, when $n = 2$. 
In \cite{Daubechies10}, it is shown that there does not exist a  compactly supported and real 
MRA wavelet when $n = 1$ and $A = 2$,  such that 
$\psi(x+1/2) = \psi(-x+1/2)$. However, in the same book,  there are examples of real and compactly supported bi-orthogonal 
wavelets $\check \psi,\check{ \tilde \psi}$ such that $ \check\psi(x+1/2) =  
\check \psi(-x+1/2)$ and 
$\check {\tilde  \psi}(x+1/2) = \check {\tilde  \psi}(-x+1/2)$.  

For general $n$ and $q = 2$, suppose  $\check \phi,  \check {\tilde \phi}$ are compactly supported and  $H$-invariant scaling functions such that the translations of $G$ form 
bi-orthogonal and  dual frames for a core space of  an MRA in $L^2(\R)$. Now $\phi, \tilde \phi \in \Xi$, 
$\langle \phi, \tilde \phi \rangle'  = 1$ and $\psi = \langle \psi , \phi \rangle' \tilde \phi = 
 \langle \psi , \tilde \phi \rangle' \phi$ for every $\psi \in C(X) \phi = C(X) \tilde \phi$.
We have that $1 = \langle \phi, \tilde \phi \rangle' = \langle \phi , \langle \tilde \phi, \tilde \phi \rangle' \phi \rangle' = \langle  \phi,  \phi \rangle' \langle \tilde \phi , \tilde \phi \rangle'$, so  $\langle \phi, \phi \rangle'$ is positive and invertible. Define $\phi_0 = \langle \phi, \phi \rangle'^{-1/2} \phi$.  We see that 
$\phi_0 \in \Xi$ satisfies $\langle \phi_0, \phi_0 \rangle' = 1$ and generates a projective multiresolution 
analysis in $\Xi$. Note that $\check \phi_0$ is not necessarily compactly supported 
if $\phi$ and $\tilde \phi$ are. We will now see that there exists an MRA wavelet in 
$U C(X) \phi_0$ with 
a symmetry as in  \cite{Daubechies10}. 

We define a  group  action $W : H \times \Xi \rightarrow \Xi$ as  
$h , \zeta \mapsto \zeta\circ h $. A computation gives that 
$\langle \zeta,\eta \rangle' \circ h = \langle W_h \zeta, W_h \eta \rangle'$.
If $W_h \phi = \phi $ and $W_h \tilde \phi = \tilde \phi$, then $W_h \phi_0 = \phi_0$. 

Following \cite{Daubechies10},  we give a formula for  an MRA wavelet 
in this case.
From now on we will use the abbreviation PMRA for a projective multiresolution analysis.

\begin{proposition} \label{q2symmetry}
  Suppose $\phi \in \Xi$, $\langle \phi, \phi \rangle' = 1$,  $q = 2$ and   $U ^j C(X) \phi$,  $j \in \Z$ form 
  a PMRA.  Let $  y_1$ be an element in the nontrivial coset of $A^T G^\bot$ in $G^\bot$ and 
  define $$\psi =  \langle  f_1 U \phi, \phi \rangle' U \phi  - \langle  U \phi , \phi \rangle' f_1 U  \phi .      $$
  Now $\check \psi$ is an MRA wavelet with mutually orthonormal translates, 
  and $$ \check \psi(h^T x - (A^T)\- y_1)   =  \check \psi(x - (A^T)\- y_1)$$ for every $h \in H$. 
\begin{proof}
  We compute 
  \begin{align*}
    & \langle \psi, \psi \rangle'  = 
      \langle U \phi,  \phi \rangle' \langle f_1 U \phi , f_1 U \phi \rangle' \langle 
     \phi , U  \phi \rangle'  \\ &  + \langle  
    f_1 U \phi  , \phi \rangle' \langle U \phi , U \phi \rangle' \langle \phi,  
    	f_1 U \phi  \rangle' 
 \\ 
     = & \langle \phi, U \phi \rangle' \langle U \phi , \phi \rangle' 
     + \langle \phi , f_1 U \phi \rangle' \langle f_1 U \phi , \phi \rangle' = \langle \phi, \phi \rangle'
  \end{align*}
and
\begin{align*}
   \langle \psi, \phi \rangle' & = - \langle U \phi, \phi \rangle' \langle f_1 U \phi , \phi \rangle'
  +  \langle f_1 U \phi , \phi \rangle' \langle U \phi , \phi \rangle' = 0.
\end{align*}

Now, $\phi$ and $\psi$ form two $C(X)$-linearly independent elements in the free $C(X)$-module 
$U C(X)\phi$ of rank two, so they generate the module.  
Finally, since $f_1 \circ h \overline f_1$ is $G$-periodic, we get   
\begin{align*}
  W_h \overline f_1 \psi & =  \overline f_1 \circ h \- \langle W_h f_1 U \phi, W_h \phi \rangle'  U \phi 
  -   \langle W_h U \phi, W_h \phi \rangle' W_h f_1 U \phi \\
  & =   \overline f_1 \circ h  \langle (f_1 \circ h \overline f_1) f_1,\phi \rangle'U \phi 
  - \langle U \phi , \phi \rangle' (f_1 \circ h \overline f_1) f_1 U \phi  = \overline f_1 \psi.
\end{align*}
If we apply the inverse Fourier transform on both sides of this equation, we obtain
$$ \check \psi(h^T x - (A^T)\- y_1)   =  \check \psi(x - (A^T)\- y_1).$$
\end{proof}
\end{proposition}

\begin{definition} \label{symdef} Let $A$ be a dilation matrix and suppose $H$ is a finite and affiliated group in $GL_n(\Z)$.   
Suppose $y_1, \dots, y_q$ form a system of representatives for the cosets of $A^T G^\bot$ in $G^\bot$ and
$y_1 \in A^T G^\bot$. 
We will say that an MRA wavelet-family $\{\check \psi_{i,j}\}_{ 1 \leq i \leq d, 2 \leq j \leq q}$ is symmetric if   
$$
	\check \psi_{i,j} ( h^T x - A^{T-1} y_j) = \check \psi_{i,j} ( x - A^{T-1} y_j) 
$$
for every $h \in H$, $1 \leq i \leq d$ and $2 \leq j \leq q$.
\end{definition}

A formula as in the statement of Proposition  \ref{q2symmetry} does not exist for $q > 2$, but we can at least give  some sufficient 
conditions for the existence of  a symmetric and orthonormal wavelet family. This is what we aim at in the next two sections.

\section{Equivalent embeddings of equivariant vector bundles}
In this section we will assume that $G$ is a finite group. 
A $G$-space is a 
topological space $X$ with a continuous group action $G \times  X  \mapsto X$. 
A continuous map between two $G$-spaces that commutes with the group 
actions is called an equivariant map.    

We will need  the notion of an equivariant simplicial complex. 
This is a  finite set, $K$,   of simplices such that 
\begin{enumerate}
\item If $s \in K$,  then the faces  of $s$ are contained in $K$, 
\item The intersection of two simplices in $K$  is either empty 
  or a face for both, 
\item $G$ acts on $K$ by simplicial maps, i.e. it 
  takes vertices of a simplex to the vertices of another simplex in $K$, 
\item 
  For every subgroup $H \subset G$  and $[v_0, \dots, v_n] \in K$, where 
  $[h_0v_0, \dots,h_n v_n] \in K$  and $h_j\in H$, there exists an $h \in H$
  such that $h v_j = h_j v_j$ for every $0 \leq j \leq n$,
\item The vertices $ \{ v_0, \dots, v_n\}$  of a simplex  in $K$ can be ordered such that  the isotropy groups satisfy  $G_{v_n} \subset \dots \subset  G_{v_0}$.

\end{enumerate}
We let $|K|$ denote the geometrical realization of the equivariant simplicial complex $K$, i.e. the union of the simplices 
embedded into $\R^m$ for a suitable $m\in \N$. When equipped with the subspace topology,  $|K|$ is a topological $G$-space with the action inherited from $K$. 
In the main theorem of  \cite{Illman}, Illman shows that if $G$ is a finite group acting smoothly on 
a smooth manifold $X$, there exists an equivariant simplicial complex $K$ and an equivariant homeomorphism 
$|K| \rightarrow X$.

Let $X$ be  a $G$-space and let $\rho : \xi \rightarrow X$ be  a vector bundle with a $G$-action. 
We say that $\xi$ is an equivariant vector bundle if  $g\rho = \rho g$ and $g | _{\xi_x}$ is a linear map onto $\xi_{gx}$ 
for every $g \in G$. 
The aim of this section is  to give conditions for when two embeddings of equivariant vector bundles are equivalent.  
Our main tool is the computation 
in \cite[4.5.3]{Hatcher} that shows that the Stiefel manifold $V_r(\C^s)$  of non-zero and orthogonal $r$-tuples in $\C^s$, 
is $2s - 2r$-connected. To us, this means that every  continuous map  $S^{n-1} \rightarrow V_r(\C^s)$
extends continuously to a map $D^n \rightarrow V_r(\C^s)$, when $n-1 \leq 2s - 2r $.

\begin{proposition} \label{cancellation}
  Suppose $X$ is an $n$-dimensional smooth manifold and
  $G$ is a finite group that acts smoothly on $X$.
  Let  $\xi$ be a  trivial $d$-dimensional $G$-equivariant vector bundle over $X$
   such that whenever  $x \in X$ and $H \subset G$ is the isotropy group of $x$,  
  then the action of $H$ on $\xi|_x$ is trivial. 
  Moreover, let  $\theta =   X \times \C^r $ be the $r$-dimensional trivial  bundle 
  over $X$ with the $G$-action defined by $g : (x,v) \mapsto (gx, v)$ for $g \in G$.
  
  If 
  $u_0,u_1 : \theta \rightarrow \xi$ are injective and equivariant 
 bundle maps and      $d \geq n/2 + r $, then there exists
  an equivariant bundle 
  automorphism  $s : \xi \rightarrow \xi$ such that $u_0 = s u_1.   $

  \begin{proof}
  	Since there exists a $G$-equivariant simplicial complex $K$ and an equivariant homeomorphism from $|K|$ to $X$, 
  		we can  assume that $X$ is a geometrical realization of an  equivariant simplicial complex.  	
  	Let $I = [0,1]$ and let $\xi^I$ denote the vector bundle with base space $X \times I$, total space  		
  	 $\xi \times I$ and projection map $( v,s) \mapsto (p(v), s)$ where $p : \xi \rightarrow X$ is the projection  from $\xi$ onto $X$.  
    Let $\zeta \subset \xi^I |_{X \times \{0,1\}}$ denote the trivial sub-bundle
    that has a copy of $u_0( \theta)$ over $X \times \{0\}$ and a copy of $u_1(\theta)$ over $X \times \{1\}$.
    If  $\tilde \zeta \subset \xi^I$ is an equivariant and trivial sub-bundle of  $\xi^I$ 
    such that $\tilde \zeta | _{X \times \{0,1\}} = \zeta$, then  $\xi^I / \tilde \zeta$ is a vector bundle over $X \times I$ such that
    	$(\xi^I / \tilde \zeta)| _{X \times \{ i\}} \simeq \xi/ u_i( \theta)$ for $i = 0,1$.
    Lemma \cite[1.6.4]{Atiyah} states that whenever $Y$ is a $G$-space and $\eta$ is an equivariant vector bundle  over 
    $Y \times I$,  then $\eta |_{Y \times \{0\}} \simeq \eta | _{Y \times \{1\}}$. This implies that 
    there exists an equivariant isomorphism
    $
    \xi/ u_0(\theta) \rightarrow \xi /u_1(\theta).
    $
     Proposition  \cite[1.6]{Atiyah} says  that short exact sequences of equivariant vector bundles split, so  
    we obtain $\xi \simeq u_i(\theta) \oplus (\xi/ u_i(\theta))$ and hence also  an equivariant bundle automorphism, 
    $s$ that intertwines $u_0$ and $u_1$.

    Pick an inner product on $\xi$.  
     The $G$-average of an inner product on $\xi$ is still an  inner product on $\xi$, 
     so we can assume that the  inner product on $\xi$ is $G$-invariant.    
  		Let $\langle \cdot , \cdot \rangle_x$ denote the restriction of our inner product on 
  			$\xi$ to the vector space $\xi|_x$. 
  			We equip $\xi^I$ with the $G$-invariant inner product that coincides with $\langle \cdot , \cdot \rangle$ on   
  			the restriction $\xi^I |_{(x,s)}$ for every $(x,s) \in  X \times I$.
       Let $s_1, \dots s_r : X \times \{0,1\} \rightarrow \zeta$
    be orthogonal and $G$-invariant sections.
    If we can extend these sections to linearly independent and  $G$-invariant sections 
    $\tilde s_1, \dots,  \tilde s_r : X \times I \rightarrow \xi ^ I$, then we  obtain a sub-bundle $\tilde \zeta \subset \xi^I$ as needed.    
    Using induction on the dimension of $X$, we will now prove that such sections do exist.

    Suppose $X$ is $0$-dimensional,  let $x \in X$ and let $H$ denote the 
    isotropy group of $x$, i.e. $H = \{ h \in G | h x = x\}$.
		Let $V_r( \xi |_x   )$ denote the Stiefel manifold of 
      orthogonal $r$-tuples in $\xi|_x$.
        The sections $s_1, \dots , s_r$, restricted to $\{x \} \times \{0,1\}$, 
      form a map $\{x\} \times \{0,1\} \rightarrow  V_r(\xi|_x)$. The space $V_r(\xi_x)$ is path connected, so 
      this   map extends continuously to a map $\{x \} \times I \rightarrow V_r(\xi_x)$.
      We obtain $r$ orthogonal sections $s_1' , \dots, s_r'  : \{x\} \times I \rightarrow  \xi^I|_{\{x\} \times I}$ 
       that 
	coincide with $s_1, \dots, s_r$ on $\{x\} \times \{0,1\}$.
	The section $s_1', \dots , s_r'    $ can be extended to orthogonal $G$-invariant sections 
        $\tilde s_1,\dots, \tilde s_r :  G x \times I
	\rightarrow  \xi ^I $. 
        Finally, since $X$ is a disjoint union of 
	such orbits, we can extend 
	$\tilde s_1, \dots, \tilde s_r$ to a family of orthogonal and $G$-invariant sections 
		$\tilde s_1, \dots , \tilde s_r : X \times I \rightarrow \xi^I$
	that coincide with $s_1, \dots , s_r$ on $X \times \{0,1\}$.
	
  Suppose $2(d-r) \geq n > 0$ and  $X$ is $n$-dimensional and let $X^{(n-1)} \subset X $ be the geometrical realization  of the 
  $G$-simplicial complex of the $(n-1)$-dimensional simplices in $X$. 
  Our induction hypothesis is that there exist $r$ orthogonal and $G$-invariant sections from 
  	  $X^{(n-1)} \times I $ to $\xi^I$
  that coincides with $s_1, \dots, s_r$ on $X^{(n-1)} \times \{0,1\}$. 
  Since they coincide on $X^{(n-1)} \times \{0, 1\}$,  we obtain $r$ orthogonal and $G$-invariant 
  	sections $\tilde  s_1, \dots, \tilde s_r : (X^{(n-1)}\times I )\cup ( X \times \{0 , 1\}) \rightarrow 
  \xi ^I$.
  Let $\sigma$ be a simplex in $X$ with isotropy group $H$.
  	The isotropy group $H$ acts trivially on the trivial bundle  $\xi ^ I | _{\sigma \times I} \simeq \sigma \times I \times \C^d  $. 
     The sections $\tilde s_1 |_{\partial \sigma \times I}, \dots, \tilde s_r|_{\partial \sigma \times I}$ 
     form a continuous map $\partial (\sigma \times I)  \rightarrow V_r( \C^d)$. 
     The space $V_r( \C^d)$ is $2(d-r)$-connected,  so this map extends continuously to a 
     map $\sigma \times I \rightarrow V_r( \C^{\dim F})$. This new map gives us  
     extensions of 
     our  sections $\tilde s_1, \dots, \tilde s_r$ to  $H$-invariant and orthogonal 
     sections in $\xi ^I | _{\sigma \times I}$. Moreover,  they  extend uniquely to $G$-invariant 
     and orthogonal sections 
     over $G \sigma$. 
     
          We can follow this procedure for the orbit of every $n$-dimensional simplex in $X$. 
     If  $\sigma$ and $\sigma'$ are two such  $n$-simplices, then 
     the chosen families of sections coincide on $\partial (\sigma \times I   ) \cap \partial ( \sigma' \times I)$, so
     we  obtain $r$  
     $G$-invariant and linearly independent sections in $\xi ^I $ that coincide with $s_1, \dots, s_r$ on $X \times \{0,1\}$.

\end{proof}
\end{proposition}

\section{Symmetries when $dq > 2$}
In this section we will apply the results from the previous section to 
give sufficient conditions for the existence of symmetric MRA wavelets when $dq > 2$.
Now $A$ is supposed to be a dilation and $H$ is an affiliated group as in definition \ref{affdef}.
$H$ acts smoothly on the smooth manifold $\T^n$, 
so Proposition \ref{cancellation} applies in this situation.  
Recall that $Q$ is the orthogonal projection onto $V_0$, $U$ is the unitary operator on $L^2(\R^n)$ defined 
as $U f = \sqrt q\-  f \circ A\- $ and $W_h$ is defined as $W_h f = f \circ h$ for $f \in \Xi$.
The main theorem in this paper gives us conditions for when there exist symmetric MRA-wavelets in the sense of Definition \ref{symdef}.  

\begin{theorem} \label{symwavelets}
Suppose $\phi_1, \dots, \phi_d \in \Xi$ are mutually orthonormal and the submodule they generate, say $V_0$, 
form a projective multiresolution analysis together with $U$. Suppose  $W_h \phi_j = \phi_j$
for every $h \in H$ and $1 \leq j \leq d$. 
Let $y_1, \dots, y_q$ be a system of representatives for the cosets of $A^T G^\bot$ in $G^\bot$ such that 
$y_1 \in A^T G^\bot$. If 
$q \geq \frac n{2d} + 1$ 
,  there exists 
an orthonormal  MRA wavelet family $\check \psi_{k,l}$ for $ 2 \leq k \leq q$ and $1 \leq l \leq d$,  
such that $\psi_{k,l} \in (1-Q)U  V_0$ and 
$$
\check \psi_{k,l}(h^T x + (A^{T})\-y_k  ) = \check \psi_{k,l}( x + (A^T)\- y_k  )
$$
for every $x \in \R^n$,  $h \in H$, $2 \leq k \leq q$ and $1 \leq l \leq d$. 
\begin{proof}
Recall the definition $f_j = e^{2 \pi i \langle A\- \cdot , y_j \rangle}$.
We assume that $y_1 = 0$, i.e. $f_1 = 1$. 
Since $f_k U \phi_l$ for $1 \leq k \leq q$ and $1 \leq l \leq d$ form a $C(X)$-basis for $U V_0 $, the map 
$S : V_0 \rightarrow U V_0$, defined by $S \sum_{j =1}^d a_j  \phi_j = \sum_{j = 1}^d a_j U  \phi_j$ 
is a well defined $C(X)$-linear map.  
Moreover, $S$ is injective and a simple computation shows that 
$SW_h = W_h S$ for every $h \in H$. We want to find an equivariant $C(X)$-linear automorphism on $UV_0$ that 
extends $S$. We identify $UV_0$ with the sections in the product bundle $X \times \C^{dq}$, so that 
$f_k U\phi_l$ is  identified with the constant section $x \mapsto (x, e_{k,l})$. 
Since $U \phi_j$ is $W_h$-invariant, we see that 
$
W_h ( f_k U \phi_l) =f_k \circ h U  \phi_l = ( f_k \circ h \overline f_k ) f_k U \phi_l, 
$
so the action of $H$ on $UV_0$ corresponds to the 
following action on $X \times \C^{dq}$:
$$
(h, ([x]_G, \sum_{k,l} \alpha_{k,l} e_{k,l})) \mapsto ( [h x']_G  ,  \sum_{k,l} \alpha_{k,l} 
 e^{2 \pi i \langle  (h-1)x' ,   y_k\rangle }e_{k,l})
$$
for an $x'$ such that $[x']_G = x  $.
Note that if $(h-1) x \in G$, then $\langle (h-1)x, y_j \rangle \in \Z$, so the isotropy group 
of $x$ acts trivially on the fibre over $x$.

Both $S$ and the inclusion $V_0 \subset UV_0$ correspond to two equivariant bundle monomorphisms 
$u_i : \theta^d \rightarrow X \times \C^{dq}$ for  $i = 1, 2$. 
Let $u_1$ be the one that corresponds to $S$. By the definition  of $S$,  there exists an  equivariant and trivial
sub-bundle in $X \times \C^{dq}$, such that the direct sum of this and $u_1(\theta^d)$ equals $X\times \C^{dq}$.

Proposition  \ref{cancellation} now applies,  so we obtain an equivariant 
bundle automorphism on $X \times \C^{dq}$ that intertwines $u_1$ and $u_2$. This bundle automorphism now corresponds to an 
equivariant $C(X)$-linear automorphism $\tilde S : U V_0 \rightarrow U V_0$, such that 
$\tilde S |_{V_0} = S$.

Equip $\oplus^{d(q-1)} C(X)$ with the ordinary $C(X)$-valued inner product. We define 
$$M : \oplus^{d(q-1)} C(X) \rightarrow (1-Q)U V_0$$ by 
$Me_{k,l} = \tilde S f_k U  \phi_l$ for 
 $2 \leq k \leq q$ and 
 $1 \leq l \leq d$. 
A computation shows that $M$ is adjointable, 
$M^* \zeta = \sum_{k,l} \langle \zeta , \tilde S f_k U \phi_l   \rangle' e_{k,l}$
 and 
$MM^* \zeta = \sum_{k,l} \langle \zeta, \tilde S f_k U \phi_l \rangle' \tilde S f_k U  \phi_l$. 
Moreover 
\begin{align*}
  W_h MM^* W_h \- \zeta & =  
  \sum_{k,l} \langle W_h\- \zeta, \tilde S f_k U \phi_l \rangle' \circ h  W_h  \tilde S f_k U \phi_l 
  \\ & =  \sum_{k,l} \langle \zeta, W_h \tilde S f_k U \phi_l \rangle'  e^{2 \pi i \langle A\- x, h^T y_k - y_k   \rangle} 
  \tilde S f_k U \phi_l \\
  & = \sum_{k,l} \langle \zeta, e^{2 \pi i \langle A\- x, h^T y_k - y_k   \rangle}       \tilde S f_k U \phi_l \rangle'  e^{2 \pi i \langle A\- x, h^T y_k - y_k   \rangle} 
  \tilde S f_k U \phi_l \\
  & = MM^* \zeta.
\end{align*}
By the proof of Proposition \ref{ortbase},  $MM^*$ is a selfadjoint and invertible operator on the 
Hilbert module $(1-Q)UV_0$. We obtain another bounded operator $(MM^*)^{-1/2}$ and
we let $\psi_{k,l} = (MM^*)^{-1/2} \tilde S f_k U \phi_l$ for 
$2 \leq k \leq q$ and $1 \leq l \leq d$. These  form an  orthonormal $C(X)$-basis for $(1-Q)UV_0$.  
The operator $(MM^*)^{-1/2}$ is a  limit of polynomials in $MM^*$ with respect to the operator norm, so we have the relation 
$W_h (MM^*)^{-1/2} = (MM^*)^{-1/2}W_h$ for every $h \in H$. 
Finally, we compute 
\begin{align*}
  W_h \overline f_k \psi_{k,l} & = \overline f \circ h (MM^*)^{-1/2}\tilde S W_h f_k U \phi_l \\
  & = \overline f_k \circ h (MM^*)^{-1/2} \tilde S e^{2\pi i \langle A\- \cdot, h^T y_k - y_k \rangle} f_k U \phi_l \\
  & = \overline f_k \circ h  e^{2\pi i \langle A\- \cdot, h^T y_k - y_k \rangle} \psi_{k,l} \\
  & = \overline f_k \psi_{k,l}.
\end{align*}
If we apply the inverse Fourier transform on both sides in this equation, we obtain  
$$
\check \psi_{k,l}(h^T x + (A^{T})\-y_k  ) = \check \psi_{k,l}( x + (A^T)\- y_k  ).
$$	 

\end{proof}
\end{theorem}

\begin{remark}
Note that since $m$ is not necessarily a Laurent polynomial, $\check \phi$ is most likely not compactly supported. In order to obtain existence 
of compactly supported and symmetric wavelets, we will probably need to consider bi-orthogonal wavelets instead of orthonormal. The problem is to give a sharp condition on $d,q,n$ such that  if the  $\phi_1,   \dots, \phi_d , \tilde \phi_1, , \dots , \tilde \phi_d \in L^2(\R^n)$  have compact supports and if the following three conditions are satisfied,  
\begin{itemize}
\item  $\tilde \phi_i \circ h^T = \tilde \phi_i$ and $\phi_i \circ h^T = \phi_i$ for every $1 \leq i \leq d$ and $h \in H$, 
\item  $
\langle  \phi_i( \cdot - k), \tilde \phi_j ( \cdot - l)\rangle = \delta_{i,j} \delta_{k,l}  
$
 for every $k,l \in \Z^n$ and $1 \leq i, j \leq d$, 

\item $W_0 = \overline{ span \{ \phi_i ( \cdot - k) | 1 \leq i \leq d, k \in \Z^n \}}$ satisfies  
\begin{itemize}
	\item 	$W_0 \subset \check U W_0$, 
	\item $\cup_{k \in \Z}  \check U^k W_0 $ is dense in $L^2(\R^n)$, 
	\item $\cap_{k \in \Z}  \check U^k  W_0  = \{0\}$. 
\end{itemize}
\end{itemize}
then there exist $\psi_{i,j} \in L^2( \R^n)$ and $\tilde \psi_{i,j} \in L^2( \R^n)$  for $1\leq i \leq d$ $1 \leq j \leq q-1$ such that: 
\begin{enumerate}
	\item The vectors   $\psi_{i,j} \in L^2( \R^n)$ and $\tilde \psi_{i,j} \in L^2( \R^n)$  for $1\leq i \leq d$ and $1 \leq j \leq q-1$
	      have compact supports. 
	\item For every $1 \leq i, i'  \leq d$ $1 \leq   j, j' \leq q-1$ and $k, l \in \Z^n$, we have $$
\langle  \psi_{i,j}( \cdot - k), \tilde \psi_{i',j'} ( \cdot - l)\rangle = \delta_{i,i'} \delta_{j, j'}  \delta_{k,l}. $$
\item The space $\check U W_0 \ominus W_0 $ is the closed linear span of 
$$      
\{   \psi_{i,j}(\cdot - k ) | 1 \leq i \leq d, 1 \leq  j \leq q - 1  , k \in \Z^n  \} 
$$
\item If $y_1, \dots , y_{q-1}$ are  representatives for the nontrivial co-sets of $\Z^n /  A^T \Z^n$, then we have the symmetry relations:
\begin{align*}
 \psi_{i,j}(h^T x + (A^{T})\-y_j  ) & =  \psi_{i,j}( x + (A^T)\- y_j  ) \\
{\tilde \psi}_{i,j}(h^T x + (A^{T})\-y_j  ) & = {\tilde  \psi}_{i,j}( x + (A^T)\- y_j  )
\end{align*}
for every $x \in \R^n$,  $h \in H$, $1 \leq j \leq q-1$ and $1 \leq i \leq d$. 
\end{enumerate}
This is a more algebraic problem than that of Theorem \ref{symwavelets}. In fact if we had a sharp algebraic version of Proposition \ref{cancellation} with free modules over the Laurent polynomials and $\C$-linear group actions, we could follow the proof of  Theorem \ref{symwavelets}  to prove existence of compactly supported and symmetric bi-orthogonal wavelets. This will be the subject of a later paper. 
\end{remark}

\section{$H$-invariant scaling functions}

 The proof of Theorem \ref{symwavelets} depends on  the nonconstructive proof of Proposition \ref{cancellation}.
 As we have seen, the  existence comes from a computation of some homotopy groups for the Stiefel manifolds.
 The pessimist would probably interpret this result as just an  absence of obstructions 
to the existence of  symmetric wavelets in a  projective multiresolution analysis. 
A more optimistic way to interpret this result is that it gives us 
 a  clue for where to  search for nice  scaling functions. 

We wish to construct a family of scaling functions for a PMRA such that  Theorem \ref{symwavelets} applies.
First we define a continuous map $r : X \rightarrow X$ with the equation 
$r \circ p(x)= p(Ax)$ for every $x \in \R^n$. Let $m \in  M_d( Lip_1(X))$. This map is often refered to as the low-pass filter.
The transfer operator associated to $m$ is a linear operator $R : M_d( Lip_1(X)) \rightarrow M_d( Lip_1(X))$ that is defined as follows:
$$
Ru (x) = q\- \sum_{y \in r\-(x)} m^*(y)u(y)m(y).
$$

We need some basic facts about how to  construct scaling functions  from such transfer operators. The following theorem is a summary of several more general results in  \cite{DuRo1}. 
\begin{theorem}[\cite{DuRo1}]
Let $w\in \C^d$ be a unitvector such $m(0) w = \sqrt q w$ and suppose the following properties are satisfied: 
\begin{enumerate}
	\item $R 1 = 1$. \label{p1}
	\item $sp(R) \cap \T =  \{1\}$. \label{p2}
	\item The geometric multiplicity of $R$'s eigenvalue $1$ equals $1$. \label{p3}
	\item $sp(q^{-1/2} m(0)) \cap \T = \{1\}$.  \label{p4}
	\item The  algebraic multiplicity of $m(0)$'s eigenvalue $\sqrt q$ equals $1$.\label{p5}
\end{enumerate}

Then  
$$
\P(x) = \lim_k q^{-k/2} m(p(A^{-1} x) ) \dots m( p(A^{-k}x)) 
$$
converges uniformly on compact sets,   the equation 
\begin{align} \label{scalingdef}
\phi_i(x)  = \langle \P(x) e_i , w \rangle
\end{align}
defines an orthonormal family $\phi_1, \dots , \phi_d \in \Xi$  and the  submodule 
$$
V_0 = C(X) \phi_1 + \dots + C(X) \phi_d \subset \Xi
$$
yields a projective multiresolution analysis $(U, \{U^kV_0\}_{k \in \Z})$in $\Xi$. 
\end{theorem}
Note that if $H$ is an affiliated group to the dilation $A$ and if $m\circ h = m $ for every $h \in H$, then 
$\P(hx) = \P(x)$ for every $h \in H$ and $x \in \R^n$. This immediatly implies that $W_h \phi_i = \phi_i$ for every $h \in H$, i.e.
the scaling functions   
$\phi_1, \dots, \phi_d$ are $H$-invariant. 

Now let $d = 1$, $n = 2$, 
$A = \begin{pmatrix} 0 & 2 \\ -2 & 0 \end{pmatrix}$ and  let $H\subset  GL_2(\Z)$ denote the group generated by $h = -1$.  
Moreover, let 
$$
m'(x) = \frac{\sqrt 2}{4} \sum_{j = 1}^4 \cos( 2 \pi \langle v_j, x \rangle)
$$
where $v_1 =  \begin{bmatrix} 1 \\ 0  \end{bmatrix}$, 
$v_2 =  \begin{bmatrix} 0  \\ 1  \end{bmatrix}$, 
$v_3 =  \begin{bmatrix} 1 \\ -1 \end{bmatrix}$, 
and  
$v_4 =  \begin{bmatrix} 2 \\ 0  \end{bmatrix}$.

We define a linear operator operator $R': C(X) \rightarrow C(X)$  with the following equation: 
$$
R'u (x) = q\- \sum_{y \in r\-(x)} \overline{m'}(y)u(y)m'(y).
$$
Let 
$F$ be a finite subset of $G^\bot = \Z^2$ and define the following set of Laurent polynomials: 
$$
\K_F =  \{ \sum_{l \in F} a_l e^{2 \pi \langle l, x \rangle } |  a_l \in \C    \}.     
$$
A short computation gives that 
$$
R' \K_F \subset \K_{(A^T)\- (F + V) \cap \Z^2},  
$$
where $V = \{ \pm v_i \pm  v_j | 1 \leq i,j \leq 4 \}$, so if 
$$J = \{k \in \Z^2 | ~\| k \| \leq \frac {  \|(A^T)\- \| \max_{l \in V    } \| l \|    }{1 - \| (A^T)\- \|      }  + 1   \},$$
then $R' \K_J \subset \K_J$. Moreover, there  exists a $K$ such that $R'^k \K_F \subset \K_J$ for every $k \geq K$.
A numerical computation with MATLAB shows that $ 1$ is the unique eigenvalue for $R' |_{\K_J   }$ in $\T$.  We can show numerically that there exist a  $u \in \K_J$ and a $c > 0$ 
such that $u(x) \geq c$ for every $x \in X$ and $u$ spans the fixed space for  $R' |_{\K_J   }$.

Suppose that $\lambda  \in \T$ and $v \in  C(X) $ such that $R v = \lambda v$. If $w$ is a Laurent polynomial such that 
$\| v - w \|  <\epsilon$, there exists a $K$ such that $k \geq K$ implies that $R'^kw \in \K_J$, i.e. 
\begin{align}
	\inf_{y \in \K_J} \|  v -  y\| & \leq \| v - \lambda^{-k} R'^k w \| \leq \| \lambda^k v - R'^k w \| \leq \| R'^k (v - w) \| \leq \sup_j \| R'^j\| \| v - w \| \\  & \leq 
	 \sup_j \| R'^j\| \epsilon.  \label{eps1}
\end{align}
Since there exists a $c >0$ such that $u(x) \geq c$ for every $x \in X$, we see as in  \cite{DuRo1}  that    $\sup_k \|R^k\|  < \infty $. 
The inequality \ref{eps1} is satisfied for an arbitrary $\epsilon > 0 $, so  $v \in \K_J$, i.e. $1$ is the unique eigenvalue for $R'$ in $\T$  and the fixed space of  $R'$ is exactly the  $1$-dimensional subspace in $\K_J$ that is fixed by $R'|_{\K_J}$.

We define $m \in Lip_1(X)$ with the following equation: 
$$
	m(x)  = u^{-1/2}(r x)m'(x) u^{1/2}(x).
$$
As noted in \cite{DuRo1}, whenever $v \in C(X)$, $\lim_k k\- \sum_{j= 1}^k  R^j v$ converges uniformly to a fixed point for $R'$ for every $v \in C(X)$. Since the fixed space of $R'$ is $1$-dimensional, we see that $u = u(0) \lim_k \sum_{j= 1}^k k\- R^j u$. Moreover, since $m' \circ h = m'$ for every $h \in H$, we obtain  $u^{-1/2} \circ h = u^{-1/2}$ and $u^{1/2} \circ h = u^{1/2}$, so by the definition of $m$, we have that $m \circ h = m$ for every $h \in H$.

We see that $m$ satisfies the properties \ref{p1}, \ref{p4} and \ref{p5}. To see that \ref{p2} and \ref{p3} are satisfied, we define  a linear automorphism $\Theta : Lip_1(X) \rightarrow Lip_1(X)$ with the equation: 
$
\Theta (v) = u^{1/2} v u^{1/2}
$. A short computation shows that $R' \Theta = \Theta R$, so the eigenspace of an eigenvalue $\lambda$ for $R'$ is isomorphic to the 
corresponding eigenspace for the operator $R$. This implies that \ref{p2} and \ref{p3} are satisfied. 
 Now \ref{scalingdef}, defines an $H$-invariant scaling function in $\Xi$ such that  Theorem 
\ref{symwavelets} applies, i.e.  this scaling function allows a symmetric MRA  wavelet family in $\Xi$.

\section{Dilations  and affiliated groups when    $n = 2$ }

Recall that a  matrix $A \in M_n(\Z)$ is said to be a dilation if 
none of  the eigenvalues of $A$ are contained in the closed unit disk and a finite subgroup $H \subset GL_n(\Z)$ is affiliated to $A$ if 
$h A = A h$ for every $h \in H$ and $H$ acts trivially on $\Z^n/A\Z^n$.
Two matrices $C_1, C_2 \in M_n(\Z)$ are equivalent if there exists an $S \in GL_n(\Z)$ such that 
$SC_1 = C_2S$. If $A$ is a dilation, $H$ is an affiliated group, and 
$S \in GL_n(\Z)$, then $SAS\-$ is another dilation with an affiliated group 
$SHS\-$. We have already seen that if $n  =  1$ and $H$ is nontrivial,  
then $A = \pm 2$ and $H \simeq \Z/2\Z$ is generated by $h = -1$.
For the case $n = 2$, we will compute all possible  pairs,  up to this equivalence. 
This list will give us a hint for where to look for examples of symmetric wavelets in $L^2(\R^2)$. 

The following result will be of great importance for us.
\begin{proposition} \cite[IX.14]{Newman}
The finite subgroups of $SL_2(\Z)$ are cyclic of order $1, 2 , 3 , 4, 6$
and are generated by a matrix  that is equivalent to one of the following:
$$
\pm 
\begin{pmatrix}
1  & 0  \\ 0 & 1
\end{pmatrix}
,
\begin{pmatrix}
0 & 1 \\  -1 & 0\\
\end{pmatrix}
,
\pm
\begin{pmatrix}
0 & 1 \\ -1&  -1  
\end{pmatrix}.
$$
\end{proposition}
\begin{lemma}
If  $h \in GL_2(\Z)$
has  order $n$, then $h$ is diagonalizable and both eigenvalues of 
$h$ are $n$'th roots of $1$. Moreover, if $\text{det}(h) = -1$ then 
$h$ is equivalent to one  and only one of the following matrices:
$$
\begin{pmatrix}
  1 & 0 \\ 0 & -1
\end{pmatrix}
,
\begin{pmatrix}
  0 & 1 \\ 1 & 0
\end{pmatrix}.
$$
\begin{proof}
If $h^n = 1$ and $\lambda$ is an eigenvalue of $h$ then 
$h^n v = \lambda^n v = v$ for an eigenvector $v$. Moreover, if $h$ is not diagonalizable then  $h$
has Jordan form  
$$
\begin{pmatrix}
  \lambda_1 & 1 \\
  0 & \lambda_2 \\
\end{pmatrix}
$$
but this matrix has not finite order. 
Suppose $\det(h)  = -1$. 
If $h$ has a complex eigenvalue, then its complex conjugate 
is also an eigenvalue for $h$. Now  
$\det(A) = |\lambda|^2 = 1 \neq -1$, i.e. 
$h$ has eigenvalues $1, -1$.

So $h$ has the form 
$$
h = 
\begin{pmatrix}
 h_{11} & h_{12} \\
 h_{21} & -h_{11} \\
\end{pmatrix}.
$$
Moreover,  $h$ has an eigenvector for the eigenvalue $1$ with integer coefficients which can be chosen to be relatively  prime, say 
$\begin{pmatrix} x \\ y \end{pmatrix}$. 
The extended Euclidean  algorithm gives us an algorithm for how to 
express the greatest common divisor of two numbers as a $\Z$-linear combination of them, i.e.
there exist two  integers $\tilde x, \tilde y$ such that 
$x\tilde x - y \tilde y = 1$. The pair $\tilde x, \tilde y$ is also known as the Bezout coefficients of $x,y$. Define 
$$
S = \begin{pmatrix} x & \tilde x \\ y & \tilde y \\ \end{pmatrix} \in  \text{GL}_2(\Z).
$$
We can compute that 
$$S\-hS = 
\begin{pmatrix}
1 & m \\
0 & -1 \\
\end{pmatrix} = h_m
$$ for an integer $m$. 
If $m$, $n$ are integers such that $m-n = 0 \mod 2$ and 
$$
T = \begin{pmatrix}
1 & \frac{m-n}{2} \\
0 & 1
\end{pmatrix},
$$
then 
$
T h_m = h_n T
$, i.e. 
$h_m$ and $h_n$ are equivalent. 
This implies that there are at most two equivalence classes of integer matrices with $\det h = -1$ 
and finite order. 
We can compute that the matrices that intertwine 
$\begin{pmatrix}
1 & 0 \\ 0 & -1
\end{pmatrix}$
and 
$\begin{pmatrix}
0 & 1\\
1 & 0
\end{pmatrix}$
are on the form 
$
\begin{pmatrix}
m & n\\
m &-n
\end{pmatrix}$.
These matrices have determinant $-2mn$, i.e. none of them are elements in $GL_2(\Z)$,
so 
there are exactly two equivalence classes of integer matrices with eigenvalues $1, -1$. 
\end{proof}
\end{lemma}

\begin{lemma}
If $h \in GL_2(\Z)$,  $\det (h) = -1$ and  $h$ generates a finite subgroup of $G$ that is affiliated to a dilation $A$,  
then $A$ and $h$ are simultaneously equivalent 
to one of the following pairs 
\begin{enumerate}
\item \label{firstp}
$$
A =  \begin{pmatrix}
   \pm 2 &   0 \\
   0  &  \pm 2 \\
 \end{pmatrix}
, 
h =  \pm \begin{pmatrix}
   1 & 0  \\
    0 &  -1 \\
 \end{pmatrix},
$$
\item \label{secondp}
$$
A =  \begin{pmatrix}
   n   &   0 \\
   0  &  \pm 2 \\
   \end{pmatrix}
, 
h = 
 \begin{pmatrix}
   1 & 0  \\
    0 &  -1 \\
 \end{pmatrix}
$$ 
with $|n| \geq 3$.

\end{enumerate} 
\begin{proof}

Suppose 
$h = 
\begin{pmatrix}
  1 & 0 \\
  0 & -1 \\
\end{pmatrix}
$. 
We can compute that if $A \in M_2(\Z)$  such that
$Ah= hA$ then,
$$
A = \begin{pmatrix}
  m & 0 \\ 
  0 & n \\
  \end{pmatrix}.
$$
Moreover, $h$ acts trivially on 
$ \Z^2 / A \Z^2 $ if and only if 
 $A\-(h-1) \in M_2(\Z)$. Now  
$$A\-(h-1) = \begin{pmatrix} 0 & 0 \\ 0 & -\frac 2 n  \end{pmatrix},$$
so $|n| = 1,2$ implies that $H$ acts trivially on $\Z^2/A\Z^2$.
Since we require that the norm of 
$A$'s eigenvalues should be strictly greater than $1$, $A$ must be one of the following matrices
$$
A = 
\begin{pmatrix}
n & 0 \\
0 & \pm 2\\
\end{pmatrix}
$$
with $|n| \geq 2$.

The equivalence imposed by the matrix $\begin{pmatrix} 0 & 1\\ 1 & 0 \end{pmatrix}$ gives that the dilations that are affiliated with 
	 $\begin{pmatrix} -1 & 0\\ 0 & 1 \end{pmatrix}$ are on the form 
	 $$
A = 
\begin{pmatrix}
\pm 2 & 0 \\
0 & n \\
\end{pmatrix}
$$
with $|n| \geq 2$.
If we let $|n| = 2$ we obtain (\ref{firstp}) and if we let $|n| \geq 3$, we obtain (\ref{secondp}).

The set of matrices that commute with  $h = \begin{pmatrix} 0 & 1 \\ 1 & 0 \end{pmatrix}$
are 
$\{\begin{pmatrix}
m & n \\ n & m
\end{pmatrix} | m,n \in \Z   \}
$. If $A$ is an affiliated  dilation  then
$$
A\-(h-I) = \det A \- 
\begin{pmatrix}
  -m+n & m -n \\
  -n +m  & n -m
\end{pmatrix} \in \text{GL}_2(\Z).
$$
Now $\det A = ( m + n ) (m -n)$ divides   $m-n$, i.e. $m+n = \pm 1$ and  $A = \begin{pmatrix} m  & m \\ m & m \end{pmatrix}  \pm \begin{pmatrix} 0 & 1 \\ 1 & 0 
\end{pmatrix}$. A computation shows that such matrices have 
either $1$ or $-1$ as an eigenvalue. 
This implies that $h$ does not commute with a proper dilation. 
\end{proof}
\end{lemma}

\begin{theorem}
If $A \in M_2(\Z)$ is a dilation and $H \subset GL_2(\Z)$ is a finite and non trivial affiliated group, 
then $A$ and $H$ are simultaneously equivalent to one of the following pairs:  

\begin{enumerate}
\item 
$A$ is one  of the following matrices
$$
 \begin{pmatrix}
   0 &  2 \\
  \pm 1 &  0 \\
 \end{pmatrix},
 \pm \begin{pmatrix}
   0  &  2 \\
  -1 &  1 \\
 \end{pmatrix},
$$
$H$ is generated by 
$$
 \begin{pmatrix}
   -1 & 0 \\ 
   0 &  -1 \\
 \end{pmatrix},
$$
$|det(A)| = 2$ and $H \simeq \Z/2\Z$.
\item
$A$ is one of the following matrices
$$
 \begin{pmatrix} 0 & 2 \\ -2 & 0  \end{pmatrix},
 \begin{pmatrix} 0 & 2 \\ -2 & \pm 2   \end{pmatrix},
 \begin{pmatrix} 0 &2  \\ 2 & \pm 2  \end{pmatrix},
 \begin{pmatrix} 2 &2  \\ 0 & -2  \end{pmatrix},
$$
$$
 \begin{pmatrix} 2 & |n|  \\ 0 & 2  \end{pmatrix},
 \begin{pmatrix} -2& |n|\\ 0 & -2  \end{pmatrix},   \text{ where }n = 0\mod 2, n \neq 0,
$$
$H$ is generated by 
$$
 \begin{pmatrix}
   -1 & 0 \\ 
   0 &  -1 \\
 \end{pmatrix},
$$
$|\det A | = 4$ and $H \simeq  \Z/2\Z$. 
\item 
$A$ is one  of the following matrices
$$
 \pm\begin{pmatrix}
   1 &  -1 \\
  1 &  1 \\
 \end{pmatrix},
$$
$H$ is generated by 
$$
 \begin{pmatrix}
   0 & 1 \\ -1&  0 \\
    \end{pmatrix}, 
$$
$det(A) = 2$ and $H \simeq \Z/4\Z$.

\item 
$A$ is one  of the following matrices
$$
 \pm\begin{pmatrix}
   2 &  1 \\
  -1 &  1 \\
 \end{pmatrix},
 \pm\begin{pmatrix}
   -1 &  -2 \\
  2 &  1 \\
 \end{pmatrix},
$$
$H$ is generated by 
$$
 \begin{pmatrix}
   -1 & -1 \\ 
   1 &  0 \\
    \end{pmatrix},
$$
$|det(A)| = 3$ and $H \simeq \Z/3\Z$.
\item 
$A$ is one  of the following matrices
$$
 \begin{pmatrix}
   \pm 2 &  0 \\
  0 &  \pm 2  \\
 \end{pmatrix},
$$
$H $ is generated by 
$$
 \pm \begin{pmatrix}
   1 & 0 \\ 
   0 &  -1 \\
    \end{pmatrix},
$$
$|det(A)| = 4$ and $H \simeq \Z/ 2 \Z \oplus \Z/2 \Z$.
\item 
$A$ is one  of the following matrices
$$
 \begin{pmatrix}
    n &  0 \\
  0 &  \pm 2  \\
 \end{pmatrix},
$$
$|n| \geq 3$, 
$H$ is generated by 
$$
  \begin{pmatrix}
   1 & 0 \\ 
   0 &  -1 \\
    \end{pmatrix},
$$
$\det(A) = \pm 2n$ and $H  \simeq \Z/2\Z$.

\end{enumerate}
\begin{proof}

Recall that we say a dilation is affiliated to $H$ if $Ah = hA$ for every $h \in H$
and $H$ acts trivially on $\Z^2/A\Z^2$. Note that this implies that $A\- (h-I) \in M_2(\Z)$ for every $h \in H$.
Suppose $h = -I$ and $A= \begin{pmatrix} a & b \\ c & d    \end{pmatrix}$ is an affiliated dilation.
Now $A\- (h-I) = \frac{ -2} { \det A} \begin{pmatrix}  d & -b \\ -c & d\end{pmatrix}$. 
If $2 \nmid \det A$, then $\det A | a, b, c ,d$, so $\det A = \det A ^2 m$ for an $m \in \Z$. 
This implies that $m ,\det A \in \{\pm 1\}$, so $A$ can not a be dilation. 
If $2 | \det A$, then  $\det A = 2 m$ for an $m \in \Z \setminus \{0\}$ such that
$m  | a, b , c ,d$. Now $\det A = m^2 n = 2m$, so $mn = 2$. This implies that $\det A = \pm 2$ or 
$\det A = \pm 4$.

We note that every  dilation $A$, such that $\det A = \pm 2$, is affiliated to $h$.  
    A classification by Lagarias and Wang \cite[Lemma 5.2]{LagariasWang} gives that every dilation with $\det A = \pm 2$ is
    equivalent to one of the following matrices: 
$$
    \begin{pmatrix}
      0 & 2 \\ 1 & 0
    \end{pmatrix}
,    \begin{pmatrix}
      0 & 2 \\ -1 & 0
    \end{pmatrix}
 ,   \pm \begin{pmatrix}
      1 & 1 \\ -1 & 1
    \end{pmatrix}
 ,   \pm \begin{pmatrix}
       0 & 2 \\ -1 & 1
    \end{pmatrix}.
     $$

A complete list of dilations $A$,  with  $\det A = \pm 4$, is given in 
\cite{KiratLau}. The dilations that are affiliated to $h$ are similar to one of the following matrices
$$
 \begin{pmatrix} 0 & 2 \\ -2 & 0  \end{pmatrix},
 \begin{pmatrix} 0 & 2 \\ -2 & \pm 2   \end{pmatrix},
 \begin{pmatrix} 0 &2  \\ 2 & \pm 2  \end{pmatrix},
 \begin{pmatrix} 2 &0  \\0  & -2   \end{pmatrix},
 \begin{pmatrix} 2 &2  \\ 0 & -2  \end{pmatrix},
$$
$$
 \begin{pmatrix} 2 & |n|  \\ 0 & 2  \end{pmatrix},
 \begin{pmatrix} -2& |n|\\ 0 & -2  \end{pmatrix}  \text{ where }n = 0 \mod 2 .
$$

If $|H| = 3$, then $H$ is  conjugate to the group generated by 
$\begin{pmatrix}0 & 1 \\ -1 & -1\end{pmatrix}$. A computation gives that 
the set of $A \in M_2(\Z)$ such that $hA = Ah$ is
$$ \{
\begin{pmatrix}
  a+b & a \\
  -a & b
\end{pmatrix} | 
~ a,b \in \Z\}.
$$
Moreover, if $A$ is an affiliated dilation,  then 
$$
A\-(h-I) = \frac{1}{\det(A)} 
\begin{pmatrix}
  a-b & 2a+b \\
  -2a-b & -a-2b
\end{pmatrix}
\in M_2(\Z)
$$
i.e. 
$\det(A) | 3b, 3a$. 
If $3 \nmid \det A$, then $\det A | a,b$. Now  $\det A = \det A ^2 m$ for an 
$m \in \Z\setminus \{0\}$, so $\det A = \pm 1$. If $3 | \det A$, there exists an
$m \in \Z$ s.t. $m | a,b,c,d$ and  $3 m = \det A = m^2n$ for an $n \in \Z \setminus \{0\}$. 
This implies that  $\det A =     \pm 3 $ or $\det A = \pm 9$.

After  some analysis,  we see that the affiliated dilations are parametrized 
by the following values of $a$ and $b$:
$$
(a,b) = \pm ( 1, 1), (1,-2),(2,-1).
$$
These values give  the following matrices 
$$
\pm \begin{pmatrix}
  2 & 1 \\
  -1 & 1
\end{pmatrix}
, 
\pm \begin{pmatrix}
  -1 & -2 \\
  2 & 1
\end{pmatrix}
,
\pm \begin{pmatrix}
  1 & -1 \\
  1 & 2
\end{pmatrix}.
$$
Moreover, if 
$S = \begin{pmatrix} 0  & 1\\ 1 & 0 \end{pmatrix}$ then 
$
S \begin{pmatrix}
  2 & 1 \\ -1 &  1 \\
\end{pmatrix} S 
=  \begin{pmatrix}
  1 & -1 \\ 1 &  2 \\
\end{pmatrix} 
$
and $
ShS = h^2$,  so we end up with the announced matrices
in the case  $|H| = 3$.

Suppose $H$ is a cyclic group of order $4$ and $H$ is generated by 
$h = \begin{pmatrix} 0 & 1 \\ -1 & 0 \end{pmatrix}$.
We compute that the set of integer matrices that commute with $h$ 
are 
$$
\{\begin{pmatrix}
  a & b \\ -b & a
\end{pmatrix} |  a,b \in \Z \}.
$$
If $A$ is an affiliated dilation then 
$$
A\-(h-I) = \frac{1}{\det(A)} 
\begin{pmatrix}
  -a+b & a+b \\ -a-b & b-a
\end{pmatrix}
\in M_2(\Z),
$$
i.e.
$\det A | 2 a , 2 b$.  
This implies that $\det A = \pm 2$ or $\det A \pm 4$.

The affiliated dilations  are parametrized by 
$(a,b) \in \{ \pm (1 ,-1   ), \pm ( 1 , 1)\} $
and  we obtain the matrices 
$$
\pm \begin{pmatrix} 1 & -1 \\ 1 & 1 \end{pmatrix}
,\pm \begin{pmatrix} 1 & 1 \\ -1 & 1 \end{pmatrix}.
$$
Note that 
$
S \begin{pmatrix} 1 & -1 \\ 1 & 1 \end{pmatrix}S 
=  \begin{pmatrix} 1 & 1 \\ -1 & 1 \end{pmatrix}
$, so we end up with the announced matrices when $H \simeq \Z/4\Z$.

Suppose $H$
is a cyclic group of order $6$ and is 
generated by 
$
h = \begin{pmatrix}
  1 & 1 \\ -1 & 0
\end{pmatrix}
$. 
If $A$ is an  affiliated dilation, it must also be affiliated to  the subgroups generated by  $h^3 = -1$ and $h^2 =  \begin{pmatrix} 0 & 1 \\ -1 & -1 \end{pmatrix}$. 
The above classification implies  that this is impossible.  
\end{proof}
\end{theorem}

\begin{acknowledgments}
The author is grateful to the referee for many valuable comments which have been incorporated into the manuscript and is also pleased to acknowledge helpful discussions with  Ola Bratteli, Tore A. Kro, Dorin E. Dutkay, Halvard Fausk, 
  Sergey Neshveyev, Nadia S. Larsen and John Rognes. 
\end{acknowledgments}

\bibliographystyle{alpha}
\bibliography{7022}

\end{document}